\documentclass[12pt]{amsart}

\newtheorem*{thm}{Theorem}
\theoremstyle{remark}
\newtheorem*{rems}{Remark}
\newtheorem{example}{Example}

\def\F{{\bf F}}
\def\P{{\bf P}}
\def\mod{\mathop{\rm mod}\nolimits}

\begin{document}
\title[Unboundedness of rational points]{Unboundedness of the number
  of rational points on curves over function fields} 

\author{Ricardo Concei\c c\~ao}
\address{Department of Mathematics, Oxford College of Emory
  University, Oxford, GA 30054, USA}
\email{rconcei@emory.edu}

\author{Douglas Ulmer}
\address{School of Mathematics, Georgia Institute of
    Technology, Atlanta, GA 30332, USA}
\email{ulmer@math.gatech.edu}

\author{Jos\'e Felipe Voloch}
\address{Department of Mathematics, University of Texas,
    Austin, TX 78712, USA}
\email{voloch@math.utexas.edu}

\begin{abstract}
We give examples of sequences of smooth non-isotrivial curves for every 
genus at least two, defined over a rational function field of positive 
characteristic,
such that the (finite) number of rational points of the curves in the sequence 
cannot be uniformly bounded.
\end{abstract}

\maketitle

The question of whether there is a uniform bound for the number of
rational points on curves of fixed genus greater than one over a fixed
number field has been considered by several authors.  In particular,
Caporaso et al. \cite{CaporasoHarrisMazur97} showed that this would
follow from the Bombieri-Lang conjecture that the set of rational
points on a variety of general type over a number field is not Zariski
dense.  Abramovich and the third author \cite{AbramovichVoloch96}
extended this result to get other uniform boundedness consequences of
the Bombieri-Lang conjecture and gave some counterexamples for
function fields. These counterexamples are singular curves that
``change genus''. They behave like positive genus curves (and, in
particular, have finitely many rational points), but are
parametrizable over an inseparable extension of the ground field.  In
\cite{AbramovichVoloch96} it is shown that, for this class of
equation, uniform boundedness does not hold. Specifically, one gets a
one-parameter family of equations which, for suitable choice of the
parameter, have a finite but arbitrarily large number of solutions.
However, a negative answer to the original uniform boundedness
question for smooth curves of genus at least two remained open in the
function field case. In this paper we provide counterexamples to this
uniform boundedness, extending constructions of the first two authors
\cite{ConceicaoThesis, UlmerLegendre} for elliptic curves.

\begin{thm}
  Let $p>3$ be a prime number and let $r$ be an odd number coprime to
  $p$.  The number of rational points over $\F_p(t)$ of the curve
  $X_a$ with equation $y^2 = x(x^r+1)(x^r+a^r)$ is unbounded as $a$
  varies in $\F_p(t)\setminus\F_p$.
\end{thm}

\begin{proof} 
  We first note that if $d=p^n+1$ and $a=t^d$, then we have a rational
  point $(x,y)=(t,t^{(r+1)/2}(t^r+1)^{d/2})$ on $X_a$.  Second, if $m$
  divides $n$ and $n/m$ is odd, then $d'=p^m+1$ divides $d$.  Setting
  $e=d/d'$, we have another rational point
  $(x,y)=(t^e,t^{e(r+1)/2}(t^{re}+1)^{d'/2})$ on $X_a$.  Thus if we
  take $n$ to be odd with many factors, we have many points.
\end{proof}

The curves given by the theorem are non-isotrivial of odd genus $r$.
To obtain counterexamples of even genus, one may proceed as follows: Let
$Y_a$ be the quotient of $X_a$ by the fixed-point-free involution
$(x,y) \mapsto (a/x,-a^{(r+1)/2} y/x^{r+1})$.  Then $X_a$ is an
unramified cover of $Y_a$, a curve of genus $(r+1)/2$. This shows that
for all $g>1$ and all but finitely many $p>2$ (the exceptions
depending on $g$), unboundedness of the number of rational points over
$\F_p(t)$ holds for curves of genus g.

One can obtain other explicit examples by a slight modification of the
argument of the theorem.  In what follows $d=p^n+1$, $m|n$ with $n/m$
odd, $d'=p^m+1$, $e=d/d'$, and $a=t^d$, as in the proof of the
theorem.

\begin{example}
  Let $p\equiv 2\mod 9$ and $n$ be divisible by 3 and a product of
  primes $\equiv 1\mod 6$. Then the curve $y^6 = x(x+1)(x+a)$ contains
  the point $(t^e,t^{e}(t^e+1)^{d'/2})$.
 \end{example}

\begin{example}
  Let $f(x)\in\F_p[x]$ be a polynomial of degree $2b$ with distinct
  roots, none of them zero.  Then the curve $y^2=f(x)x^{2b}f(a/x)$ has
  the point $(t^e,t^{be}f(t^e)^{d'/2})$.
\end{example}

\begin{example}
  Let $r$ be a prime satisfying $p\equiv r-1\mod r^2$. Let $n$ be
  divisible by $r$ and primes $\equiv 1 \mod r(r-1)$.  Then
  $(t^e,t^{2e/r}(t^e+1)^{d'/r})$ is a point on the curve
  $y^r=x(x+1)(x+a)$.  This curve has a simple Jacobian, since the ring
  of integers of the $r$-th cyclotomic field acts as endomorphisms of
  the Jacobian (see \cite{Zarhin06}, Theorem 3.1).
\end{example}

\begin{rems}\mbox{}
  The curve in Example 1 and the curve in the theorem cover the
  Legendre elliptic curve $E$ studied in \cite{UlmerLegendre}. Thus
  their Jacobians have $E$ as a factor, and consequently have
  unbounded Mordell-Weil rank as $a$ varies. On the other hand, $E$ is
  not a factor of the Jacobian of the curve $Y_a$, nor of the
  Jacobians of the curves in Examples 2 and 3.  Nonetheless, by the
  main result of \cite{BuiumVoloch96}, the rank of the Mordell-Weil
  group of the Jacobian of these curves is also unbounded as $a$
  varies.
\end{rems}

Let $X$ be the smooth projective surface with affine model
$$y^2=x(x^r+1)(x^r+t^{r}).$$  
By \cite{CaporasoHarrisMazur97} the fibration $X \to \P^1$, $(x,y,t)
\mapsto t$ has a fibered power which covers a variety of general type.
However, since this fibration is defined over a finite field, this
variety of general type will also be defined over a finite field and
can have a Zariski dense set of $\F_p(t)$-rational points, so the rest
of the argument of \cite{CaporasoHarrisMazur97} does not apply.
\bigskip

\emph{ Acknowledgments:} We would like to thank Bjorn Poonen for his
suggestions, in particular for suggesting the involution that leads to
the curves $Y_a$.

\bibliography{database}{}
\bibliographystyle{alpha}

\end{document}